\documentclass{article}
\usepackage{amsmath}[1996/11/01]
\usepackage{amssymb,amsthm,amsxtra}

\def\[#1\]{\begin{equation}#1\end{equation}}
\makeatletter
\def\beq{%
   \relax\ifmmode
      \@badmath
   \else
      \ifvmode
         \nointerlineskip
         \makebox[.6\linewidth]%
      \fi
      $$
   \fi
}
\def\eeq{%
   \relax\ifmmode
      \ifinner
         \@badmath
      \else
         $$
      \fi
   \else
      \@badmath
   \fi
   \ignorespaces
}

\def\enddisplaymath{\eeq\global\@ignoretrue}
\makeatother

\newtheorem{thm}{Theorem}
\newtheorem{cor}[thm]{Corollary}
\newtheorem{lem}[thm]{Lemma}

\newtheorem{conj}[thm]{Conjecture}
\theoremstyle{remark}
\newtheorem{eg}{Example}
\newtheorem*{rem}{Remark}
\newtheorem{rems}{Remark}[thm]

\numberwithin{equation}{section}
\numberwithin{thm}{section}

\DeclareMathOperator{\Exp}{{\bf E}}
\DeclareMathOperator{\Var}{Var}
\newcommand{\Z}{\mathbb Z}
\newcommand{\R}{\mathbb R}
\newcommand{\N}{\mathbb N}

\begin{document}
\title{A mean identity for longest increasing subsequence problems}
\author{Eric Rains\footnote{AT\&T Research, New Jersey, 
rains@research.att.com}}
\date{October 16, 2000}
\maketitle
\begin{abstract}
We show that a wide variety of generalized increasing subsequence problems
admit a one parameter family of extensions for which we can exactly compute
the mean length of the longest increasing subsequence.  By the nature of
the extension, this gives upper bounds on the mean in the unextended model,
which turn out to be asymptotically tight for all of the models that have
so far been analyzed.  A heuristic analysis based on this fact gives not
just the asymptotic mean but also the asymptotic scale factor, again
agreeing with all known cases.
\end{abstract}

\section{Introduction}

In \cite{SpohnP2}, Pr\"ahofer and Spohn consider a certain polynuclear
growth (PNG) model with stationary initial conditions, and show that it
maps to the following increasing subsequence problem:

Let $t$ be a positive real number.  Pick a random set of points in the unit
square $[0,1]\times [0,1]$ as follows.  On the left and bottom edges, take
a Poisson process of mean $t$; inside the square, take a Poisson process of
mean $t^2$.  (Thus our total mean is $t^2+2t$.)  A sequence of these points
is ``increasing'' if we have $x\le x'$, $y\le y'$ whenever $(x,y)$ and
$(x',y')$ are consecutive points in the sequence; the length of the
sequence is defined to be the number of points.  The problem is then to
determine the asymptotic distribution of the length of the longest
increasing subsequence.  (Note that without the extra points on the left
and bottom edge, this is just the standard Poisson model for increasing
subsequences of random permutations \cite{AD}.)

Pr\"ahofer and Spohn then observe \cite{SpohnP3} that the stationarity of
the initial conditions can be used to show that the length of the longest
increasing subsequence has mean exactly $2t$.  This fact is striking for
two reasons.  The first is that the mean in the standard model is rather
complicated; it is thus surprising that a fairly simple extension gives
rise to an exact formula for the extended mean.  The second is that since
adding points can only help the longest increasing subsequence, we conclude
that $2t$ is an upper bound on the mean in the standard model.  This bound
is quite tight; indeed, in the standard model, the mean takes the form
$2t+O(t^{1/3})$.  (We could also derive this upper bound from the (strictly
stronger) result of \cite{Pilpel} that the expected length of the longest
increasing subsequence of a permutation of length $n$ is at most
$2\sqrt{n}$; the present method is more generally applicable, however.)

The first object of the present paper is to generalize this fact in a
number of ways.  It turns out, for instance, that the standard Poisson
model admits a one-parameter family of extensions with explicit means; this
further extends to give explicit information about the moment generating
function in a neighborhood of this family.  Also, we can replace the
Poisson model by the generalized model considered in \cite{BR1} (based in
turn on a model of Johansson \cite{kurtj:shape}).

Our other object is to explore the asymptotic relations between the
extended models and the unextended models.  It turns out that by a careful
(if heuristic) analysis, we can use the moment generating function
identities to determine not just the asymptotic mean of the longest
increasing subsequence length, but also the asymptotic scale factor.  This
gives a uniform prescription for the scaling information, agreeing with the
results of all of the cases that have so far been analyzed.

Section 1 defines the models of interest, as well as a certain continuous
limiting case.  Section 2 gives a short, algebraic proof of the moment
generating function and mean identities; this is followed by a somewhat
more complicated, but also more enlightening combinatorial proof in Section
3.  Finally, Section 4 considers the asymptotic consequences, giving
explicit conjectures for the asymptotics of the general case.

{\bf Acknowledgements}  The author would like to thank J. Baik,
Y. Baryshnikov, and W. Whitt for helpful conversations.
\section{Extended growth models}

The model we will be considering is a generalization of a model considered
by Johansson \cite{kurtj:shape}, combining the generalizations of
\cite{BR1}, \cite{SpohnP2}, and \cite{BR:png}.  We define a ``parameter
set'' $p$ to be a triple $(t,q,r)$, where $t$ is a nonnegative number, and
$q$ and $r$ are sequences of nonnegative numbers with
\[
\sum_i q_i + \sum_i r_i < \infty;
\]
we will write such a parameter set as $t{:}q/r$, omitting $t$ if $t=0$ and
omitting any trailing 0's from $q$ and $r$.  From the existence of the sum,
we conclude that
\[
Q(p):=\sup_i(q_i)\quad\text{and}\quad R(p):=\sup_i(r_i)
\]
are well-defined and attained.  Given a parameter set $p$, we define two
functions
\begin{align}
H(z;p) &= e^{tz} \prod_i (1-z q_i)^{-1}(1+z r_i)\\
E(z;p) &= e^{tz} \prod_i (1+z q_i)(1-z r_i)^{-1};
\end{align}
$H(z;p)$ converges for $|z|<Q(p)^{-1}$, while $E(z;p)$ converges for
$|z|<R(p)^{-1}$.  We will say that two parameter sets $p_+=t^+{:}q^+/r^+$
and $p_-=t^-{:}q^-/r^-$ are {\it compatible} if $Q(p_+)Q(p_-)<1$,
$R(p_+)R(p_-)<1$.  Note that we then have
\[
H(p_+;p_-):=
e^{t^+t^-}
\prod_i e^{t^- (q^+_i + r^+_i)}
e^{t^+ (q^-_i + r^-_i)}
\prod_{i,j}
(1-q^+_i q^-_j)^{-1}(1-r^+_i r^-_j)^{-1}(1+q^+_ir^-_j)(1+r^+_iq^-_j)
<\infty.
\]

Let $(\Z^+)'$ be a disjoint copy of $\Z^+$, and consider the set
$\Omega:=[0,1]\uplus \Z^+\uplus (\Z^+)'$.  Then we associate to a pair $p_+$,
$p_-$ of compatible parameter sets a random multiset $M(p_+,p_-)\subset
\Omega\times\Omega$ as follows.  We let $P(t)$ denote a Poisson random
variable of parameter $t$, $g(q)$ denote a geometric random variable with
parameter $q$, and $b(q)$ denote a random variable which is 0 with
probability ${q\over 1+q}$ and 1 with probability ${1\over 1+q}$.
\begin{itemize}
\item{} On $[0,1]\times [0,1]$, we choose $P(t^+t^-)$ i.i.d. uniform points,
\item{} On $[0,1]\times i$, we choose $P(t^+q^-_i)$ i.i.d. uniform points,
\item{} On $[0,1]\times i'$, we choose $P(t^+r^-_i)$ i.i.d. uniform points,
\item{} On $i\times j$, we have multiplicity $g(q^+_iq^-_j)$,
\item{} On $i\times j'$, we have multiplicity $b(q^+_ir^-_j)$,
\item{} On $i'\times j'$, we have multiplicity $g(r^+_ir^-_j)$,
\end{itemize}
and so on, with all of the multiplicities chosen independently.

Choose a pair of total orderings (denoted $<_+$ and $<_-$) on $\Omega$
compatible with the usual ordering on $[0,1]$.  A subsequence of a multiset
$M$ in $\Omega\times\Omega$ (a sequence of points $(x_i,y_i)$ from $M$ with
no point occuring more often than its multiplicity) is ``increasing'' if we
always have $x_i\le_+ x_{i+1}$, $y_i\le_- y_{i+1}$, subject to the further
conditions
\begin{align}
x_i=x_{i+1} &\implies x_i\notin \Z^+\\
y_i=y_{i+1} &\implies y_i\notin \Z^+
\end{align}
In other words, the sequence must be strictly increasing along rows and
columns from $\Z^+$.  We then define a sequence $\lambda_i(M)$ by setting
\[
\sum_{1\le i\le l} \lambda_i(M)
\]
equal to the size of the longest subsequence of $M$ which is a union of $l$
increasing subsequences.  We then define
$\lambda(p_+,p_-):=\lambda(M(p_+,p_-))$ (with the latter notation preferred
when particular multisets are being compared).  We recall the following
results from \cite{BR1}:

\begin{thm}\label{thm:reorder}
For any pair $p_+$, $p_-$ of compatible parameter sets, $\lambda(p_+,p_-)$
is a random partition, finite with probability 1.  The distribution of
$\lambda(p_+,p_-)$ is independent of the choice of total orderings on
$\Omega$.
\end{thm}

\begin{thm}
For any pair $p_+$, $p_-$ of compatible parameter sets, we have the identity
\[
\Pr(\lambda_1(p_+,p_-)\le l)
=
H(p_+;p_-)^{-1}
\Exp_{U\in U(l)} \det(H(U;p_+)H(U^\dagger;p_-)).
\]
\end{thm}

\begin{rem}
As in \cite{BR1}, this is a formal integral, defined by analytic
continuation from the region $Q(p_\pm)<1$.
\end{rem}

For our purposes, we need to extend the model slightly, by adding a special
row and column to the random multiset.  Extend $\Omega$ to $\Omega^+$ by
adding a new element, denoted $\Sigma$, and extend the total orderings
so that $\Sigma$ is smallest in both orderings.  Then we define a new
random multiset $M(p_+,p_-;\alpha_+,\alpha_-)$ as follows.
\begin{itemize}
\item{} On $\Omega\times\Omega$, we take $M(p_+,p_-)$,
\item{} On $[0,1]\times \{\Sigma\}$, we choose $P(\alpha_- t^+)$
i.i.d. uniform points,
\item{} On $(i,\Sigma)$, we have multiplicity $b(\alpha_- q^+_i)$,
\item{} On $(i',\Sigma)$, we have multiplicity $q(\alpha_- r^+_i)$,
\end{itemize}
and so on, and allow increasing subsequences to be weakly increasing in the
new row and column.  Note that the point $(\Sigma,\Sigma)$ has multiplicity
fixed at 0; otherwise, the new model would simply be a special case of the
old model.  By the argument in \cite{BR:png}, we have:

\begin{thm}\label{thm:integral}
For any compatible pair $p_+$, $p_-$ of parameter sets and any $\alpha_+$,
$\alpha_-$ with $\alpha_+ R(p_-)<1$ and $\alpha_- R(p_+)<1$, 
we have
\begin{multline}
\Pr(\lambda_1(p_+,p_-;\alpha_+,\alpha_-)\le l)
=\\
E(\alpha_+;p_-)^{-1}
E(\alpha_-;p_+)^{-1}
H(p_+;p_-)^{-1}
[D_l(p_+,p_-;\alpha_+,\alpha_-)-\alpha_+\alpha_-D_{l-1}(p_+,p_-;\alpha_+,\alpha_-)],
\end{multline}
where
\[
D_l(p_+,p_-;\alpha_+,\alpha_-)
=
\Exp_{U\in U(l)} \det((1+\alpha_+ U)(1+\alpha_-
U^\dagger)H(U;p_+)H(U^\dagger;p_-)).
\]
\end{thm}

Aside from intrinsic interest (\cite{BR:png}), this new model in principle can
give us some information about the original model $\lambda_1(p_+,p_-)$,
since we have the coupling
\[
\lambda_1(M(p_+,p_-))\le \lambda_1(M(p_+,p_-;\alpha_+,\alpha_-)).
\]
Thus it is particularly interesting that, as we shall see, there is an
exact formula for
\[
\Exp(\lambda_1(p_+,p_-;\alpha,\alpha^{-1}))
\]
whenever $R(p_+)<\alpha<R(p_-)^{-1}$.

We will also consider a continuous limiting case of the above model,
combining the exponential limit of \cite{kurtj:shape} and the heavy-traffic
limit of queueing theory \cite{GlynnWhitt} (studied in the present context
in \cite{TracyWidom:words} and \cite{Baryshnikov}).

The parameters for the continuous model consist of a pair of sequences
$\rho^\pm_i\in \R\cup\{\infty\}$, a nonnegative real $u$, and two numbers
$a_\pm\in \R\cup\{\infty\}$, subject to the convergence constraint
\[
\sum_i {1\over \rho^\pm_i -z} < \infty,\ \text{for}\ z<\inf(\rho^\pm),
\]
and the compatibility constraints
\[
\inf(\rho^+)+\inf(\rho^-)>0,\quad\quad
\inf(\rho^+)+a_->0,\quad\quad
a_+ + \inf(\rho^-)>0.
\]
We omit any trailing $\infty$'s from $\rho^\pm$.  Define an infinite matrix
$\Lambda_{ij}$, $0\le i,j$ by
\[
\Lambda_{ij} = \cases
0 & (i,j)=(0,0)\\
{1\over \rho^+_i+\rho^-_j},&(i,j)\ne (0,0),
\endcases
\]
with the convention that $\rho^\pm_0 = \alpha_\pm$,
and let $M$ be a random matrix filled with independent exponential random
variables, such that $M_{ij}$ has mean $\Lambda_{ij}$.  Also, for $0\le
i$ such that $\rho^+_i<\infty$, let $B_i$ be a Brownian motion on $[0,1]$
with
\[
\Exp(B_i(x)) = -u \rho^+_i x\quad\quad \Var(B_i(x)) = u x;
\]
by convention, $\rho^\pm_0 = a_\pm$.  Then we define a random sequence
$\chi(\rho^+,u{:}\rho^-;a_+,a_-)$ in terms of increasing paths in
$\N\times (\N\uplus [0,1])$, where the contribution of an interval
$i\times[x_i,y_i]$ is $B_i(y_i)-B_i(x_i)$ and the contribution of a
point $i\times j$ is $M_{ij}$.  When $u>0$, we require that
the union of $k$ increasing paths used to define $\chi_k$ have total
Lebesgue measure $k$ in $\N\times [0,1]$, and that the paths contain no
point in $i\times [0,1]$ when $\rho^+_i = \infty$.  If no such path
exists, we set $\chi_k = -\infty$.

\begin{thm}
For any valid choice of $\rho^+$, $u{:}\rho^-$, $a_+$, $a_-$, of
compatible parameter sets, $\chi(\rho^+,u{:}\rho^-;a_+,a_-)$ is a random
nonincreasing sequence in $\R\cup \{-\infty\}$, and is nonnegative
if $u=0$.  The distribution is invariant under reordering of the sequences
$\rho^\pm$.
\end{thm}

\begin{proof}
That the sequence $\chi$ is nonincreasing, invariant under reordering,
and nonnegative when $u=0$ follows from the fact that it is a scaled
limit of $\lambda$ (see below); it suffices therefore to show that
$\chi_1<\infty$ with probability 1.  We first observe that
\[
\chi_1(\rho^+,u{:}\rho^-;a_+,a_-)
\le
\chi_1(\rho^+,u{:}\rho^-;a_+,\infty)
+
\chi_1(\rho^+,u{:}\rho^-;\infty,a_-).
\]
Since decreasing $a_\pm$ cannot decrease $\chi_1$, we conclude that it
suffices to prove finiteness when $-\inf(\rho^\mp)<a_\pm<\inf(\rho^\pm)$
and $a_+ + a_- > 0$.  But then Corollary \ref{cor:continuous} below
applies, expressing $\Exp(e^{(a_+ + a_-)\chi_1})$ as an infinite product.
Moreover, the conditions on $\rho^\pm$ suffice to force convergence of this
product, and thus $\Exp(e^{(a_+ + a_-)\chi_1})<\infty$.  But this
immediately implies $\chi_1<\infty$, as desired.
\end{proof}

As we alluded to above, this is a limiting case of the discrete model:
\[
\chi(\rho^+,u{:}\rho^-;a_+,a_-)
=
\lim_{t\to\infty}
{
\lambda(/e^{-\rho^+/t},ut^2{:}/e^{-\rho^-/t};e^{-a_+/t},e^{-a_-/t})
-ut^2
\over
t}
,
\]
with the limit taken in the sense of distribution.
This corresponds to the facts that if $x$ is an exponential random variable
of mean $1/m$, then $\lfloor {x\over l}\rfloor$ is a geometric random
variable of parameter $e^{-m/l}$, and that as the parameter tends
to infinity, Poisson processes converge (with proper scaling) to Brownian
motion.

The main significance of the continuous model is that it contains two
classical matrix ensembles as special cases.  Let $0_n$ denote the finite
sequence consisting of $n$ copies of $0$.  We obtain the Gaussian Unitary
Ensemble from $\chi(0_n,1{:};\infty,\infty)$ that is,
$\chi(0_n,1{:};\infty,\infty)$ is distributed as the (ordered) eigenvalues
of an $n\times n$ Hermitian Gaussian matrix (extended to an infinite
sequence by adding $-\infty$).  (This is essentially proved in
\cite{TracyWidom:words}; to be precise, they show that one obtains
traceless GUE under a constraint equivalent to $\sum_i B_i(1)=0$; it
follows easily that without this constraint, one obtains ordinary GUE.)
Similarly, we obtain the Laguerre Unitary Ensemble from
$\chi(0_{n_+},0_{n_-};\infty,\infty)$; that is, the distribution of the
singular values of a $n_+\times n_-$ complex Gaussian matrix.
\section{An algebraic proof}

Let $p_+$, $p_-$ be a pair of compatible parameter sets.

\begin{lem}
$D_l(\alpha_+,\alpha_-):=D_l(p_+,p_-;\alpha_+,\alpha_-)$ is a polynomial,
satisfying the identity
\[
D_l(\alpha_+,\alpha_-) = (\alpha_+\alpha_-)^l D_l(\alpha^{-1}_-,\alpha^{-1}_+).
\]
\end{lem}

\begin{proof}
This follows immediately from the corresponding fact for
$\det((1+\alpha_+ U)(1+\alpha_- U^\dagger))$.
\end{proof}

For $\alpha_+<R(p_-)^{-1}$, $\alpha_-<R(p_+)^{-1}$, define
$L(\alpha_+,\alpha_-) = \lambda_1(p_+,p_-;\alpha_+,\alpha_-)$;
then

\begin{thm}\label{thm:mgf}
For $R(p_+)<\alpha_+<R(p_-)^{-1}$ and $R(p_-)<\alpha_-<R(p_+)^{-1}$,
\[
\Exp((\alpha_+\alpha_-)^{-L(\alpha_+,\alpha_-)})
=
E(\alpha_+;p_-)^{-1} E(\alpha_-;p_+)^{-1}
E(\alpha^{-1}_+;p_+) E(\alpha^{-1}_-;p_-)
\]
\end{thm}

\begin{proof}
By Theorem \ref{thm:integral} above,
\begin{align}
\Pr(L(\alpha_+,\alpha_-)\le l)
&=
E(\alpha_+;p_-)^{-1} E(\alpha_-;p_+)^{-1}
H(p_+;p_-)^{-1}
[D_l(\alpha_+,\alpha_-)-\alpha_+\alpha_- D_{l-1}(\alpha_+,\alpha_-)]\\
&=
(\alpha_+\alpha_-)^l
E(\alpha_+;p_-)^{-1} E(\alpha_-;p_+)^{-1}
H(p_+;p_-)^{-1}
[D_l(\alpha^{-1}_-,\alpha^{-1}_+)-D_{l-1}(\alpha_-^{-1},\alpha_+^{-1})].
\end{align}
Then
\begin{align}
\sum_{0\le l\le k} (\alpha_+\alpha_-)^{-l} &\Pr(L(\alpha_+,\alpha_-)=l)\notag\\
&=
\sum_{0\le l\le k}
(\alpha_+\alpha_-)^{-l}
[\Pr(L(\alpha_+,\alpha_-)\le l)-\Pr(L(\alpha_+,\alpha_-)\le l-1)]\\
&=
(\alpha_+\alpha_-)^{-k-1} \Pr(L(\alpha_+,\alpha_-)\le k)
+
{\alpha_+\alpha_- -1\over \alpha_+\alpha_-}
\sum_{0\le l\le k} (\alpha_+\alpha_-)^{-l} \Pr(L(\alpha_+,\alpha_-)\le l)\\
&=
E(\alpha_+;p_-)^{-1} E(\alpha_-;p_+)^{-1}
H(p_+;p_-)^{-1}\notag\\
&\phantom{{}={}\quad}
\left[
{1\over \alpha_+\alpha_-}
\Bigl(
D_k(\alpha_-^{-1},\alpha_+^{-1})-D_{k-1}(\alpha_-^{-1},\alpha_+^{-1})
\Bigr)
+
{\alpha_+\alpha_--1\over \alpha_+\alpha_-}
D_k(\alpha_-^{-1},\alpha_+^{-1})
\right]\\
&=
E(\alpha_+;p_-)^{-1} E(\alpha_-;p_+)^{-1}
H(p_+;p_-)^{-1}
[
D_k(\alpha_-^{-1},\alpha_+^{-1})
-(\alpha_+\alpha_-)^{-1}D_{k-1}(\alpha_-^{-1},\alpha_+^{-1})
]\\
&=
E(\alpha_+;p_-)^{-1} E(\alpha_-;p_+)^{-1}
E(\alpha_+^{-1};p_+) E(\alpha_-^{-1};p_-)
\Pr(L(\alpha_-^{-1},\alpha_+^{-1})\le k),
\end{align}
where the last step is valid since $\alpha_-^{-1}<R(p_-)^{-1}$
and $\alpha_+^{-1}<R(p_+)^{-1}$.

The theorem then follows by taking the limit $k\to\infty$.
\end{proof}

Taking a limit as $\alpha_+\alpha_-\to 1$, we obtain:

\begin{cor}\label{cor:mean}
Whenever $R(p_+)<\alpha<R(p_-)^{-1}$,
\[
\Exp(L(\alpha,\alpha^{-1}))
=
{\alpha E'(\alpha;p_-)\over E(\alpha;p_-)}
+
{\alpha^{-1} E'(\alpha^{-1};p_+)\over E(\alpha^{-1};p_+)}.
\]
\end{cor}

Remark.  If $p=t{:}q/r$, then
\[
{\alpha E'(\alpha;p)\over E(\alpha;p)}
=
\alpha t
+
\sum_i {\alpha q_i\over 1+\alpha q_i}
+
\sum_i {\alpha r_i\over 1-\alpha r_i}
\]
We also observe that
\[
{\alpha E'(\alpha;p_-)\over E(\alpha;p_-)}
+
{\alpha^{-1} E'(\alpha^{-1};p_+)\over E(\alpha^{-1};p_+)}
=
\alpha {d\over d\alpha}
\log {E(\alpha;p_-)\over E(\alpha^{-1};p_+)}
\]

In the continuous limit, we write
$
X(a_+,a_-) := \chi_1(\rho^+,u{:}\rho^-;a_+,a_-).
$
Taking the appropriate limit gives:

\begin{cor}\label{cor:continuous}
Let $\rho^+$ and $u{:}\rho^-$ be compatible continuous parameter sets.  Then
whenever $-\inf(\rho^-)<a_+<\inf(\rho^+)$ and $-\inf(\rho^+)<a_-<\inf(\rho^-)$,
\[
\Exp(e^{(a_+ + a_-)X(a_+,a_-)})
=
e^{u(a_-^2 - a_+^2)/2}
\prod_i {\rho^+_i+a_-\over \rho^+_i-a_+}
\prod_i {\rho^-_i+a_+\over \rho^-_i-a_-}.
\]
Whenever $-\inf(\rho^-)<a<\inf(\rho^+)$,
\[
\Exp(X(a,-a))
=
-ua + \sum_i {1\over \rho^-_i+a} + \sum_i {1\over \rho^+_i-a}.
\]
\end{cor}
\section{A combinatorial proof}

Fix parameters as in the previous section, and set
\begin{align}
N_+  &= |M(p_+,p_-;\alpha_+,\alpha_-) \cap (\{\Sigma\}\times \Omega)|\\
N_-  &= |M(p_+,p_-;\alpha_+,\alpha_-) \cap (\Omega\times \{\Sigma\})|
\end{align}
Then we observe
\[
\Exp((\alpha_+\alpha_-)^{-N_+})
=
{E(\alpha_-^{-1};p_-)\over E(\alpha_+;p_-)}\quad\quad
\Exp((\alpha_+\alpha_-)^{-N_-})
=
{E(\alpha_+^{-1};p_+)\over E(\alpha_-;p_+)},
\]
and for $\alpha_+=\alpha_-^{-1}=\alpha$,
\[
\Exp(N_+)
=
{E'(\alpha;p_-)\over E(\alpha;p_-)}\quad\quad
\Exp(N_-)
=
{E'(\alpha^{-1};p_+)\over E(\alpha^{-1};p_+)}.
\]
So we can restate Theorem \ref{thm:mgf} and Corollary \ref{cor:mean} as
\begin{align}
\Exp((\alpha_+\alpha_-)^{-L(\alpha_+,\alpha_-)})
&=
\Exp((\alpha_+\alpha_-)^{-N_+-N_-})\\
\Exp(L(\alpha_+,\alpha_-))
&=
\Exp(N_+ +N_-)
\end{align}
We give a direct proof of this fact, for $\alpha_+\alpha_-\le 1$:

\begin{proof}
Let $\alpha'<\alpha_-$.  Then $\alpha_+\alpha'<1$, so we can extend
$M(p_+,p_-;\alpha_+,\alpha')$ by adjoining $(\Sigma,\Sigma)$ with
multiplicity $N_0$ of distribution $g(\alpha_+\alpha')$; denote the
resulting random multiset by $M'$.  But then by Theorem \ref{thm:reorder},
we can change the total ordering $<_-$ so that $\Sigma$ becomes maximal
instead of minimal.  We then find that
\[
\Sigma\times \Omega,(\Sigma,\Sigma),\Omega\times\Sigma
\]
induces an increasing subsequence of $M'$ with respect to this new ordering;
thus
\[
N_0+N_++N_-
\le
\lambda_1(M').
\]
On the other hand, this is the only maximal increasing subsequence that passes
through $(\Sigma,\Sigma)$; any other maximal increasing subsequence can have
size at most $N_++N_-+\lambda_1(M(p_+,p_-))$.  We thus find
\[
\lambda_1(M')
\le
N_++N_-
+
\max(\lambda_1(M(p_+,p_-)),N_0).
\]
But $\lambda_1(M')$ is distributed as
$N_0 + \lambda_1(M(p_+,p_-;\alpha_+,\alpha'))$ (since before the reordering
every maximal increasing subsequence passes through $(\Sigma,\Sigma)$,
and $N_0$ is independent of $\lambda_1(M(p_+,p_-;\alpha_+,\alpha'))$.
So if we take the expectations and subtract/divide by the contribution of
$N_0$, we find that we need only show
\begin{align}
\lim_{\alpha'\to (1/\alpha)^-}
\Exp(\max(\lambda_1(p_+,p_-),N_0))-\Exp(N_0)
&=
0\\
\lim_{\alpha'\to (\alpha_-)^-}
\Exp((\alpha_+\alpha_-)^{-N_0})^{-1}
\Exp((\alpha_+\alpha_-)^{-\max(\lambda_1(p_+,p_-),N_0)})
&=
1.
\end{align}

\begin{lem}
Let $X$ be a nonnegative-integer-valued random variable with finite
first moment, and let $Y$ be an independent geometric random variable
of parameter $t$.  Then
\[
\lim_{t\to 1^-} \Exp(\max(X,Y))-\Exp(Y) = 0.
\]
Similarly, for $s<1$, if $\Exp(s^{-X})$ is finite, then
\[
\lim_{t\to s^-} \Exp(s^{-Y})^{-1} \Exp(s^{-\max(X,Y)}) = 1.
\]
\end{lem}

\begin{proof}
\begin{align}
\Exp(\max(X,Y)-Y)
&=
\Exp((X-Y) \Pr(Y<X))\\
&=
\Exp(X) + {t\over 1-t} (\Exp(t^X)-1)\\
&\to 0.
\end{align}
Similarly,
\begin{align}
\Exp(s^{-\max(X,Y)})
&=
\Exp(s^{-X} \Pr(Y<X) + \Pr(Y\ge X) \Exp(s^{-Y}|Y\ge X))\\
&=
\Exp(s^{-X}) + {t(s-1)\over t-s} \Exp((t/s)^X).
\end{align}
and thus
\begin{align}
\lim_{t\to s^-}
\Exp(s^{-Y})^{-1}\Exp(s^{-\max(X,Y)})
&=
\lim_{t\to s^-}
{t-s\over s(t-1)} \Exp(s^{-X}) + {t(s-1)\over s(t-1)} \Exp((t/s)^X)\\
&=
1
\end{align}
\end{proof}

The theorem then follows from the following lemma, since
\[
\alpha_+\alpha_->R(p_+)R(p_-)
\]
by assumption.
\end{proof}

\begin{lem}
For all $z>R(p_+)R(p_-)$, $\Exp(z^{-\lambda_1(p_+,p_-)})$ is finite.
In particular, since $1>R(p_+)R(p_-)$, $\lambda_1(p_+,p_-)$ has moments
of all orders.
\end{lem}

\begin{proof}
An increasing subsequence in $M(p_+,p_-)$ can pass through a point on a strict
row or column at most once; thus $\lambda_1(M(p_+,p_-))$ is unchanged if we
remove any excess multiplicity in those rows and columns.  Let $M^o$ be
the resulting multiset, then
\[
\lambda_1(M(p_+,p_-))=\lambda_1(M^o)\le |M^o|.
\]
It will thus suffice to prove that $\Exp(z^{-|M^o|})<\infty$.
But the moment generating function of $|M^o|$ is $f(z)/f(1)$, where
\[
f(z)
=
e^{-t^+t^-z}
\prod_i e^{t^+ q^-_i z} e^{t^+ r^-_i z} e^{q^+_i t^- z} e^{r^+_i t^- z}
\prod_{i,j} (1-r^+_ir^-_j z)^{-1} (1+q^+_ir^-_j z) (1+r^+_iq^-_j z) ((1-q^+_iq^-_j)+q^+_iq^-_j z)
\]
This product converges to an analytic function with no pole inside the
open disc $|z|<(R(p_+)R(p_-))^{-1}$, and thus the result follows.
\end{proof}
\section{Asymptotic consequences}

Since $\lambda_1(p_+,p_-;\alpha_+,\alpha_-)$ is nondecreasing in $\alpha_+$
and $\alpha_-$, we obtain the following bound:

\begin{thm}
For any compatible parameters $p_+$, $p_-$,
\[
\Exp(\lambda_1(p_+,p_-))
\le
\inf_{R(p_+)<\alpha<R(p_-)^{-1}}
\Exp(\lambda_1(p_+,p_-;\alpha,\alpha^{-1})).
\label{eq:averagebound}
\]
\end{thm}

For instance, in the purely Poisson case, $p_+=p_-=t{:}/$, we find
\[
\Exp(\lambda_1(p_+,p_-))\le \inf_{\alpha>0} (\alpha+\alpha^{-1})t
= 2t.
\]
This bound is remarkably tight; indeed, we have \cite{BDJ1}
\[
\lambda_1(p_+,p_-) = 2t-O(t^{1/3+\epsilon}).
\]
This suggests the following conjecture:

\begin{conj}\label{conj:weak}
Fix parameters $p_+$, $p_-$, and define
\[
m(\alpha;p_+,p_-)
=
\Exp(\lambda_1(p_+,p_-;\alpha,\alpha^{-1})).
\]
Then
\[
\lim_{n\to \infty}
n^{-1} \lambda_1(p_+^n,p_-^n)
=
\inf_{R(p_+)<\alpha<R(p_-)^{-1}} m(\alpha;p_+,p_-),
\label{eq:subadditive}
\]
with probability 1, where
\begin{align}
(t{:}q/r)^n &:= (nt){:}q^n/r^n\\
q^n &:= q_1,q_1,\dots q_1,q_2,q_2,\dots q_2,\dots.
\end{align}
\end{conj}

\begin{rems}
Roughly speaking, this is an analogue of the law of large numbers.  As
such, it can most likely be strengthened considerably (considering
different sequences of parameter sets than just $p_+^n$, $p_-^n$).
See, for instance, the result of \cite{Seppalainen1}.
\end{rems}

\begin{rems}
The existence of the limit \eqref{eq:subadditive} follows from
superadditivity and the bound \eqref{eq:averagebound}.
\end{rems}

This has been verified in a number of special cases (see below).
In each case, we in fact find that
\[
{\lambda_1(p_+^n,p_-^n)-\mu n\over n^{1/3}}
\]
converges to a limit distribution.

Fix parameters $p_+$, $p_-$.  An increasing subsequence of
$M(p_+,p_-;\alpha,\alpha^{-1})$ cannot include points from
both $\{\Sigma\}\times\Omega$ and $\Omega\times\{\Sigma\}$.  We would thus
expect that for $\alpha$ ``large'', the typical longest increasing
subsequence will avoid $\Omega\times\{\Sigma\}$ entirely.  In particular,
we would expect
\[
\Exp(\lambda_1(p_+,p_-;\alpha,\alpha^{-1}))
\sim
\Exp(\lambda_1(p_+,p_-;\alpha,0))
\]
whenever $N_+ \gg N_-$.
For asymptotic purposes, this condition is simply that
$\alpha>\tilde\alpha_+$, where $\tilde\alpha_+$ minimizes
$m(\alpha;p_+,p_-)$. (We also define $\tilde\alpha_- = (\tilde\alpha_+)^{-1}$,
which of course minimizes $m(\alpha;p_-,p_+)$).

In particular, this tells us that $\tilde\alpha_\pm$ is a critical point;
if $\alpha_+<\tilde\alpha_+$ and $\alpha_-<\tilde\alpha_-$, we have
\[
\Exp(\lambda_1(p_+,p_-;\alpha_+,\alpha_-))
\sim
\Exp(\lambda_1(p_+,p_-;\tilde\alpha_+,\tilde\alpha_-)),
\]
while if either is greater, the mean is determined by the dominant
parameter.

This behaviour is, of course, confirmed by the analysis of \cite{BR:png}, in
which the asymptotics for general $\alpha_\pm$ are determined for the
Poisson case $p=p'=t{:}/$ and the Johansson case $p=p'=/\sqrt{q}^n$
(where $\sqrt{q}^n$ is the finite sequence consisting of $n$ copies of
$\sqrt{q}$).  In both cases, we obtain the same behaviour near the critical
point.  This suggests that for general parameters there should exist
constants $\mu$, $\sigma$, and $\sigma_\pm$ so that the following holds:

\indent
If we fix $w_+$, $w_-$, and define
\[
\alpha_\pm = \tilde{\alpha}_\pm \exp(-{2w_\pm\over \sigma_\pm n^{1/3}}),
\]
then as $n\to\infty$,
\[
{\lambda_1(p^n_+,p^n_-;\alpha_+,\alpha_-) - \mu n\over \sigma n^{1/3}}
\]
converges to the distribution $H(w_+,w_-)$ (\cite{BR:png}).

We recall the following information about the distribution $H(w_+,w_-)$:

\begin{lem}
Let $X$ be distributed as $H(w_+,w_-)$.  Then
$\Exp(\exp(2 (w_+ + w_-) X)) = \exp({8\over 3} (w_+^3 + w_-^3))$.
If $w_+ = - w_- = w$, then $\Exp(X) = 4w^2$.
\end{lem}

Only the latter equation was actually shown in \cite{BR:png}, but essentially
the same calculation gives the first equation as well.  Furthermore, in
the cases that have been fully analyzed, this is precisely the analogue
for $H$ of Theorem \ref{thm:mgf} and Corollary \ref{cor:mean} above.
This suggests that we compare the results in general.

On the one hand we compute
\[
\log(\Exp((\alpha_+\alpha_-)^{-\lambda_1(p^n_+,p^n_-;\alpha_+,\alpha_-)}))
=
2 \mu n^{2/3} ({w_+\over \sigma_+} + {w_-\over \sigma_-})
+
\log(
\Exp(
\exp(2 w_+ X \sigma/\sigma_+)
\exp(2 w_- X \sigma/\sigma_-)
)),
\]
with $X$ distributed according to $H(w_+,w_-)$ in the limit.
Thus to retain the analogy, we must have $\sigma_+=\sigma_-=\sigma$.
On the other hand, we have:
\begin{multline}
\log(\Exp((\alpha_+\alpha_-)^{-\lambda_1(p^n_+,p^n_-;\alpha_+,\alpha_-)}))
=\\
2 n^{2/3} (\theta g)(\tilde\alpha_+) {w_+ + w_-\over \sigma}
+
2 n^{1/3} (\theta^2 g)(\tilde\alpha_+) {w^2_- - w^2_+\over \sigma^2}
+
{4\over 3} (\theta^3 g)(\tilde\alpha_+) {w^3_- + w^3_+\over
\sigma^3}
+ O(n^{-1/3}).
\end{multline}
where $(\theta f)(z) = z {d\over dz} f(z)$ and $g(z) = \log E(z;p_-) - \log
E(z^{-1};p_+)$.  Thus, in particular, $(\theta g)(z) = m(z;p_+,p_-)$.
Comparing the asymptotics, we find
\[
\sigma = ((\theta^3 g)(\tilde\alpha_+)/2)^{1/3}
       = (\tilde\alpha^2_+ m''(\tilde\alpha_+;p_+,p_-)/2)^{1/3}.
\]

Similar considerations (based on part (iv) below) give us the scale factors
for $\alpha_+>\tilde\alpha_+$, thus giving us the following conjecture (see
\cite{BR:png} for the definitions of the limiting distributions):

\begin{conj}\label{conj:asympt}
Fix parameters $p_+$, $p_-$, define $\tilde\alpha_\pm$ as above, and further
define
\begin{align}
\mu          &= m(\tilde\alpha_+;p_+,p_-)&
\sigma       &= (\tilde\alpha_+^2 m''(\tilde\alpha_+;p_+,p_-)/2)^{1/3}\\
\mu_+(z) &= m(z;p_+,p_-) &
\sigma_+(z) &= (z m'(z;p_+,p_-))^{1/2},\\
\mu_-(z) &= m(z;p_-,p_+) &
\sigma_-(z) &= (z m'(z;p_-,p_+))^{1/2}.
\end{align}
Assume that $\tilde\alpha_\pm\notin \{0,\infty\}$, $\sigma>0$,
$\sigma_+(\alpha_+)^2>0$ for $\tilde\alpha_+<\alpha_+<R(p_-)^{-1}$, and
$\sigma_-(\alpha_-)^2>0$ for $\tilde\alpha_-<\alpha_-<R(p_+)^{-1}$.  Then we have
the following limiting distributions as $n\to\infty$:
\begin{itemize}
\item{(i)} If $0\le \alpha_+<\tilde\alpha_+$ and
$0\le \alpha_-<\tilde\alpha_-$ are fixed, then
\[
{\lambda_1(p_+^n,p_-^n;\alpha_+,\alpha_-)-\mu n\over
\sigma n^{1/3}} \to F_{\text{\rm GUE}}
\]
\end{itemize}
Near the critical point, set $w_\pm$ by $\alpha_\pm = \tilde\alpha_\pm
\exp(-2w_\pm/\sigma n^{1/3})$.
\begin{itemize}
\item{(ii)} If $w_\pm$ and $0\le \alpha_\mp<\tilde\alpha_\mp$ are
fixed,
\[
{\lambda_1(p_+^n,p_-^n;\alpha_+,\alpha_-)-\mu n\over
\sigma n^{1/3}}
\to G(w_\pm).
\]
\item{(iii)} If $w_+$ and $w_-$ are fixed,
\[
{\lambda_1(p_+^n,p_-^n;\alpha_+,\alpha_-)-\mu n\over
\sigma n^{1/3}}
\to
H(w_+,w_-).
\]
\end{itemize}
Finally (the Gaussian regime), let $\alpha^0_+$ and
$\alpha^0_-$ be fixed such that $\alpha^0_\pm > \tilde\alpha_\pm$ and
$\mu_+(\alpha^0_+) =
\mu_-(\alpha^0_-)$, and set $\alpha_\pm =
\alpha^0_\pm \exp(x_\pm/\sigma_\pm(\alpha^0_\pm) n^{1/2})$.
\begin{itemize}
\item{(iv)} If $x_\pm$ and $0\le \alpha_\mp<\alpha^0_\mp$ are
fixed,
\[
{\lambda_1(p_+^n,p_-^n;\alpha_+,\alpha_-)-\mu_\pm(\alpha^0_\pm) n\over
n^{1/2}}
\to
N(x_\pm,\sigma_\pm(\alpha^0_\pm)^2).
\]
\item{(v)} If $x_+$ and $x_-$ are fixed,
\[
{\lambda_1(p_+^n,p_-^n;\alpha_+,\alpha_-)-\mu_\pm(\alpha^0_\pm) n\over
n^{1/2}}
\to
\max(N(x_+,\sigma_+(\alpha^0_+)^2),N(x_-,\sigma_-(\alpha^0_-)^2)).
\]
\end{itemize}
\end{conj}

\begin{rems}
This, in turn, is analogous to the central limit theorem, so again can
probably be strengthened considerably (although not nearly to the same
extent as Conjecture \ref{conj:weak} most likely can).  In particular,
it is presumably sufficient for the parameters $\alpha_\pm$, $w_\pm$,
$x_\pm$ to tend to limits as appropriate, rather than simply be fixed.
\end{rems}

\begin{rems}
We recall that
\[
\mu_+(\alpha) = \Exp(N_+) + \Exp(N_-)
\]
at $\alpha_+ = \alpha$, $\alpha_- = \alpha^{-1}$.  Similarly,
\[
\sigma_+(\alpha) = \Var(N_+) - \Var(N_-)
\]
\end{rems}

\begin{rems}
The analogous conjectures for models of the other symmetry types
(\cite{BR1}, \cite{BR2}, \cite{BR3}) are straightforward.  We note in
particular that when $p_+=p_-=p$, we find $\alpha m'(\alpha) = -\alpha^{-1}
m'(\alpha^{-1})$, and $\alpha m'(\alpha)>0$ whenever $\alpha>1$.  So the
hypotheses of the above conjectures hold in such cases, with
$\tilde\alpha_\pm=1$.
\end{rems}

In the continuous limit, we make a similar conjecture; the main difference
is that the model is nonincreasing in the parameters, not nondecreasing, so
the $F_{\text{\rm GUE}}$ region is now $a_\pm > \tilde{a}_\pm$.  The scale
factors are:
\begin{align}
\mu          &= m_c(\tilde{a}_+;\rho^+,u{:}\rho^-)&
\sigma       &= (m''_c(\tilde{a}_+;\rho^+,u{:}\rho^-)/2)^{1/3}\\
\mu_+(z)    &= m_c(z;\rho^+,u{:}\rho^-) &
\sigma_+(z)^2  &= -m'_c(z;\rho^+,u{:}\rho^-),\\
\mu_-(z)    &= \mu_+(-z)&
\sigma_-(z)^2 &= -\sigma_+(-z)^2
\end{align}
where we define
\[
m_c(z;\rho^+,u{:}\rho^-)
=
-uz + \sum_i {1\over \rho^-_i + a} + \sum_i {1\over \rho^+_i - a}
\]
Near the critical point, we take $a^\pm = \tilde{a}^\pm
+2w^\pm/\sigma n^{1/3}$, while in the Gaussian regime, we take
$a^\pm = a^\pm_0 - x^\pm/\sigma_\pm(a^\pm_0) n^{1/2}$.

As remarked above, Conjecture \ref{conj:asympt} was proved in \cite{BR:png}
(with the exception of parts (iv) and (v), which are straightforward using
the argument in section 7 of \cite{BR2}) for the cases $p_\pm=t{:}/$ and
$p_\pm=/\sqrt{q}^n$.  The only other known results are for the case
$\alpha_\pm=0$; the references in the following examples refer to this
case alone.

\begin{eg}
If we take $p_\pm=t{:}/$, we have
\[
\mu_\pm(z)      = (z+z^{-1}) t,\quad
\sigma_\pm(z)   = ((z-z^{-1}) t)^{1/2},\quad
\tilde\alpha    = 1,\quad
\mu             = 2t,\quad
\sigma          = t^{1/3}.
\]
This corresponds to the classical case of increasing subsequences of random
permutations, studied in \cite{BDJ1}.
\end{eg}

\begin{eg}
If we take $p_\pm=/(\sqrt{q})^{n_\pm}$, with $n_+/n_-$ tending
to a constant, we have
\[
\mu_+(z) = 
{\sqrt{q}z n_+\over 1-z\sqrt{q}}
+
{\sqrt{q}n_-\over z-\sqrt{q}},\quad\quad
\sigma_+(z)^2    =
{\sqrt{q}z n_+\over (1-z\sqrt{q})^2}
-
{\sqrt{q}z n_-\over (z-\sqrt{q})^2}
\]
and thus
\[
\tilde\alpha_+ = {\sqrt{qn_+}+\sqrt{n_-}\over \sqrt{qn_-}+\sqrt{n_+}},\quad
\mu      = {q(n_++n_-)+2\sqrt{qn_+n_-}\over 1-q},\quad
\sigma   =
{(qn_+n_-)^{1/6}
(1+\sqrt{qn_+\over n_-})^{2/3}
(1+\sqrt{qn_-\over n_+})^{2/3}
\over 1-q}.
\]
This model was analyzed in
\cite{kurtj:shape}, along with the continuous (Laguerre) limit, in which case we have
\[
\mu_+(z) = {1\over 2} [{n_+\over 1+2z} + {n_-\over 1-2z}],\quad\quad
\sigma_+(z)^2    =
{4 n_+\over (1+2z)^2} - {4 n_-\over (1-2z)^2}
\]
and thus
\[
\tilde{a}_+ = {\sqrt{n_+} - \sqrt{n_-}\over 2(\sqrt{n_+} + \sqrt{n_-})},\quad
\mu         = (\sqrt{n_+}+\sqrt{n_-})^2,\quad
\sigma      = (n_+n_-)^{-1/6} (\sqrt{n_+}+\sqrt{n_-})^{4/3}.
\]
\end{eg}

\begin{eg}
If we take $p_+=/1^n$, $p_-=t{:}/$, we have
\[
\mu_+(z) = z t + {n\over z-1},\quad\quad
\sigma_+(z)^2 = z t - {z n\over (z-1)^2},
\]
and thus
\[
\tilde\alpha_+ = 1 + \sqrt{n/t},\quad
\mu      = t + 2\sqrt{tn},\quad
\sigma   = \sqrt{t} n^{-1/6} (1+\sqrt{n/t})^{2/3}.
\]
This model, corresponding to weakly increasing subsequences of random
words, was studied in \cite{kurtj:plancherel} and \cite{TracyWidom:words}.
In the continuous limit (corresponding to the $n\times n$ GUE
\cite{TracyWidom:words}), we have:
\[
\mu    = 2\sqrt{n},\quad
\sigma = n^{-1/6}
\]
These are precisely the scale factors required to make the largest
eigenvalue of an $n\times n$ Gaussian Hermitian matrix tend to the limit
$F_{\text{GUE}}$.
\end{eg}

\begin{eg}
If we take $p_+=(\sqrt{q})^{n_+}/$, $p_-=(\sqrt{q})^{n_-}/$, with
$n_+/n_-$ tending to a constant, we have
\[
\mu_+(z)              = 
{\sqrt{q}z n_+\over 1+z\sqrt{q}}
+
{\sqrt{q}n_-\over z+\sqrt{q}}\quad\quad
\sigma_+(z)^2     =
{\sqrt{q}z n_+\over (1+z\sqrt{q})^2}
-
{\sqrt{q}z n_-\over (z+\sqrt{q})^2}.
\]
Here we have three cases.  If $qn_+\ge n_-$, then $\tilde\alpha_+ = 0$, and if
$qn_-\ge n_+$, then $\tilde\alpha_+=\infty$; in either case, the above
conjectures do not apply (indeed, in those cases one expects the limiting
distribution to be atomic, $\lambda_1(p_+,p_-)=\min(n_+,n_-)$).  Otherwise,
\[
\tilde\alpha_+ = {\sqrt{n_-}-\sqrt{qn_+}\over \sqrt{n_+}-\sqrt{qn_-}},\quad
\mu      = {2\sqrt{qn_+n_-} - q(n_++n_-)\over 1-q},\quad
\sigma   =
{(qn_+n_-)^{1/6}
(1-\sqrt{qn_+\over n_-})^{2/3}
(1-\sqrt{qn_-\over n_+})^{2/3}
\over 1-q}.
\]
The mean in this model was derived in \cite{Seppalainen2}; the refined
asymptotics of a symmetrized version was studied in \cite{Baik:vicious}.
\end{eg}

\begin{eg}
If we take $p_+=(\sqrt{q})^{n_+}/$, $p_-=/(\sqrt{q})^{n_-}$, with
$n_+/n_-$ tending to a constant, we have
\[
\mu_+(z)              = 
{\sqrt{q}z n_+\over 1-z\sqrt{q}}
+
{\sqrt{q} n_-\over z+\sqrt{q}},\quad
\sigma_+(z)^2      =
{\sqrt{q}z n_+\over (1-z\sqrt{q})^2}
-
{\sqrt{q}z n_-\over (z+\sqrt{q})^2}.
\]
There are two cases.  If $qn_+\ge n_-$, then $\tilde\alpha_+ = 0$, and the
above remark applies.  Otherwise
\[
\tilde\alpha_+ = {\sqrt{n_-}-\sqrt{qn_+}\over \sqrt{n_+}+\sqrt{qn_-}},\quad
\mu      = {2\sqrt{qn_+n_-} + q(n_--n_+)\over 1+q},\quad
\sigma   = {(qn_+n_-)^{1/6} (1-\sqrt{qn_+\over n_-})^{2/3} (1+\sqrt{qn_-\over n_+})^{2/3}
\over 1+q}.
\]
The mean in this model was first derived in \cite{Seppalainen3}; the
fluctuations have been analyzed in section 5 of \cite{kurtj:plancherel}.
\end{eg}

\begin{eg}
If we take $p_+=1^n/$, $p_-=t{:}/$, we have
\[
\mu_+(z) = z t + {n\over z+1},\quad\quad
\sigma_+(z)^2 = z t - {z n\over (z+1)^2}
\]
If $n\le t$, then $\tilde\alpha_+=0$; otherwise we have
\[
\tilde\alpha_+ = \sqrt{n/t} - 1,\quad
\mu      = 2 \sqrt{tn} - t,\quad
\sigma   = (tn)^{1/6} (1-\sqrt{t/n})^{2/3}.
\]
This corresponds to strictly increasing subsequences of random words,
and was studied to a small extent in \cite{TracyWidom:words}.
\end{eg}

\medskip
Note that the pathological case $\tilde\alpha_+=0$
(resp. $\tilde\alpha_+=\infty$) can only occur if all the rows (resp. columns)
are strict.

\end{document}